\documentclass[a4paper,11pt]{amsart}
\usepackage{cases, cite}
\usepackage{amsfonts,latexsym,rawfonts,amsmath,amssymb,amsthm}
\usepackage[plainpages=false]{hyperref}
\usepackage{graphicx}
\usepackage{setspace}

\RequirePackage{color}

\theoremstyle{plain}
  \newtheorem{theorem}{Theorem}[section]
  \newtheorem{lemma}[theorem]{Lemma}
  \newtheorem{proposition}[theorem]{Proposition}
  \newtheorem{conjecture}[theorem]{Conjecture}
  \newtheorem{corollary}[theorem]{Corollary}

\theoremstyle{remark} 
  \newtheorem*{remark}{Remark}

\theoremstyle{definition}

\newtheorem*{claim}{Claim}

\numberwithin{equation}{section}

\DeclareMathOperator*{\vol}{vol}
\DeclareMathOperator*{\intr}{int}

\begin{document}

\title[Local and global versions of the $L^p$-Brunn-Minkowski inequality]{Equivalence of the local and global versions of the $L^p$-Brunn-Minkowski inequality}

\author{Eli Putterman}
\address{School of Mathematical Sciences, Tel Aviv University, Tel Aviv, 66978, Israel.}
\email{putterman@mail.tau.ac.il}



\date{\today}



\begin{abstract}
By studying $L^p$-combinations of strongly isomorphic polytopes, we prove the equivalence of the $L^p$-Brunn-Minkowski inequality conjectured by B\"or\"oczky, Lutwak, Yang and Zhang to the local version of the inequality studied by Colesanti, Livshyts, and Marsiglietti and by Kolesnikov and Milman, settling a conjecture of the latter authors. In addition, we prove the local inequality in dimension $2$ for $p = 0$, yielding a new proof of the log-Brunn-Minkowski inequality in the plane.
\end{abstract}

\maketitle

\baselineskip=16.4pt
\parskip=3pt

\section{Introduction}

The log-Brunn-Minkowski inequality, conjectured by \cite{BLYZ} and proven there in dimension $2$, arises as one aspect of the Brunn-Minkowski-Firey theory established by Lutwak \cite{L2,L3}. This theory has seen extensive development since its inception; see, e.g., \cite{BLYZ1,CLYZ,LYZs,LYZt,LYZo, Z, Z1,Z2}. A comprehensive exposition of the theory may be found in Chapter 9 of the second edition of the book of Schneider \cite{S}.

The classical Brunn-Minkowski theory deals with the Minkowski sum of convex bodies, which can be presented in the form
\begin{equation}\label{mink_sum}
K + L = \bigcap_{u \in S^{n - 1}} \{x: \langle x, u\rangle \le h_K(u) + h_L(u)\}
\end{equation}

If $K, L$ are convex bodies in $\mathbb R^n$, their Minkowski sum satisfies the famous Brunn-Minkowski inequality
\begin{equation}\label{bm_ineq}
\vol(K + L)^{\frac{1}{n}} \ge \vol(K)^{\frac{1}{n}} + \vol(L)^{\frac{1}{n}}
\end{equation}
See \cite{Gar} for an elegant exposition of this inequality and its many applications in diverse areas of mathematics.


Firey \cite{F2} defined $L_p$-combinations of convex bodies in analogy to the classical Minkowski combinations:
\begin{equation}\label{firey_sum}
(1 - \lambda) K +_p \lambda L = \bigcap_{u \in S^{n - 1}} \{x: \langle x, u\rangle \le ((1 - \lambda) h_K(u)^p + \lambda h_L(u)^p)^{\frac{1}{p}}\}
\end{equation}

For $p \ge 1$, the support function of $(1 - \lambda) K +_p \lambda L$ is simply $((1 - \lambda) h_K(u)^p + \lambda h_L(u)^p)^{\frac{1}{p}}$, but this is not the case for $p \le 1$. 

As $p \to 0$, $(\lambda h_K(u)^p + (1 - \lambda) h_L(u)^p)^{\frac{1}{p}}$ approaches $h_K(u)^\lambda h_L(u)^{1 - \lambda}$, so the $L^0$ or log-Minkowski combination of $K$ and $L$ is naturally defined as
\begin{equation}\label{log_sum}
(1 - \lambda) K +_o \lambda L = \bigcap_{u \in S^{n - 1}} \{x: \langle x, u\rangle \le h_K(u)^\lambda h_L(u)^{1 - \lambda}\}
\end{equation}

For $p > 1$, Firey established an analogue of the Brunn-Minkowski inequality:
\begin{equation}\label{bmf_ineq}
\vol((1 - \lambda) \cdot K +_p \lambda \cdot L)^{\frac{p}{n}} \ge (1 - \lambda) \vol(K)^{\frac{p}{n}} + \lambda \vol(L)^{\frac{p}{n}}
\end{equation}


For $p < 1$, it is easily seen that an analogue of \eqref{bmf_ineq} does not hold for all pairs of convex bodies, even in dimension $1$.

However, restricting to centrally symmetric bodies, B\"or\"oczky, Lutwak, Yang and Zhang \cite{BLYZ} made the following conjecture, known as the $L^p$-Brunn-Minkowski inequality:

\begin{conjecture}\label{pBM_conj} For $p \in (0, 1)$ and any two centrally symmetric convex bodies $K, L \subset \mathbb R^n$, and any $\lambda \in [0, 1]$, 
\begin{equation}\label{pBM_ineq}
\vol((1 - \lambda) K +_p \lambda L) \ge (1 - \lambda) \vol(K)^{\frac{p}{n}} + \lambda\vol(L)^{\frac{p}{n}}
\end{equation}
\end{conjecture}

In the case $p = 0$, this becomes the following, which is called the log-Brunn-Minkowski conjecture:

\begin{conjecture}\label{logBM_conj} For any two centrally symmetric convex bodies $K, L \subset \mathbb R^n$ and any $\lambda \in [0, 1]$, 
\begin{equation}\label{logBM_ineq}
\vol((1 - \lambda) K +_o \lambda L) \ge \vol(K)^{1 - \lambda} \vol(L)^{\lambda}
\end{equation}
\end{conjecture}


B\"or\"oczky, Lutwak, Yang and Zhang \cite{BLYZ} showed that Conjecture \ref{logBM_conj} implies Conjecture \ref{pBM_conj} for any $p > 0$. Furthermore, \cite{BLYZ} established Conjecture \ref{logBM_conj} for bodies in the plane.

Several works have treated the log-Brunn-Minkowski conjecture since \cite{BLYZ}. Saroglou proved the log-Brunn-Minkowski conjecture in the case where both bodies are unconditional \cite{SC1}, and also demonstrated its equivalence to the (B)-conjecture for uniform measures \cite{SC1, SC2}. Ma \cite{M} provided an alternative proof of Conjecture \ref{logBM_conj} in the plane. Rotem \cite{R} observed that the conjecture for complex convex bodies $K_0,K_1 \subset \mathbb{C}^n$ follows from a more general theorem of Cordero--Erausquin \cite{CE}.

Recently, Colesanti, Livshyts, and Marsiglietti \cite{CLM} obtained a stability version of the log-Brunn-Minkowski conjecture near the Euclidean ball:

\begin{theorem}Let $R \in (0, \infty)$ and $\varphi \in C^2(S^{n - 1})$ be even and strictly positive. Then there exists a sufficiently small $a > 0$ such that for every $\epsilon_1 ,\epsilon_2 \in  (0 ,a)$ and for every $\lambda \in [0,1]$, 
\begin{equation}\label{local_logBM}
\vol((1 - \lambda) K_1 +_o \lambda K_2) \ge \vol(K_1)^{1 - \lambda} \vol(K_2)^{\lambda}
\end{equation}
where $K_i$ is the convex body with support function $R \varphi^{\epsilon_i}$.
\end{theorem}

This result was improved in the subsequent paper by Colesanti and Livshyts \cite{CL}. The method of proof, in brief, is to consider the function $f(\lambda) = \log \vol((1 - \lambda) B_2^n +_o \lambda K)$, where $h_K = R e^\varphi$, and show that \eqref{local_logBM} is equivalent to the inequality $f''(0) \le 0$. One then computes $f''$ explicitly, obtaining an integral inequality on the sphere which is proven using the theory of spherical harmonics.

Kolesnikov and Milman \cite{KM} have studied so-called local $L^p$-Brunn-Minkowski inequalities for $C^2_+$ bodies. These inequalities have several equivalent formulations; we cite one which is reminiscent of Minkowski's second inequality (and in fact strengthens it):

\begin{conjecture}\label{local_pBM_KM} Given a $C^2_+$ centrally symmetric convex body $K$ and $p \in [0, 1)$, the following inequality holds for all even functions $z \in C^2(S^{n - 1})$:
\begin{equation} \label{local_pBM_KM_eq}
\frac{1}{\vol(K)} V(z h_K[1], K[n - 1])^2 \geq \frac{n-1}{n-p} V(z h_K[2], K[n - 2]) + \frac{1-p}{n-p} V(z^2 h_K[1], K[n - 1])
\end{equation}
\end{conjecture}

Here, $V$ denotes the mixed volume, which in the setting of $C^2_+$ bodies can be defined for $n$-tuples of functions as well as convex bodies (see, e.g., \cite[\S 4.1]{KM}).

Kolesnikov and Milman showed that for a given $p$, the $L^p$-Brunn-Minkowski inequality (Conjecture \ref{pBM_conj}) implies Conjecture \ref{local_pBM} \cite[Lemma 3.4]{KM}, and also showed that Conjecture \ref{local_pBM} implies the $L^p$-Brunn-Minkowski inequality for pairs of bodies satisfying a certain condition \cite[Proposition 3.9]{KM}:

\begin{proposition}\label{geodesic_pBM}
Assume Conjecture \ref{local_pBM} holds for $p \in [0, 1)$, and let $K_0, K_1$ be $C^2_+$ symmetric convex bodies such that $K_\lambda = (1 - \lambda) K_0 +_p \lambda K_1$ is $C^2_+$ for all $\lambda \in [0, 1]$. Then the $L^p$-Brunn-Minkowski inequality holds for $K_0, K_1$ and any $\lambda$:
\begin{equation}
\vol(K_\lambda) \ge (1 - \lambda) \vol(K_0)^{\frac{p}{n}} + \lambda\vol(K_1)^{\frac{p}{n}}
\end{equation}
\end{proposition}

The assumption of the proposition holds in particular if $K_1$ lies in a sufficiently small neighborhood of $K_0$, which is the reason that inequalities of the form \eqref{local_pBM_eq} are termed local $L^p$-Brunn-Minkowski inequalities. Kolesnikov and Milman conjectured that for a given $p$, the local $L^p$-Brunn-Minkowski inequality is in fact equivalent to the $L^p$-Brunn-Minkowski inequality \cite[Conjecture 3.8]{KM}.

Using Riemannian geometry methods, Kolesnikov and Milman were able to prove Conjecture \ref{local_pBM_KM} for $p \in [1 - c n^{-\frac{3}{2}}, 1]$ for an absolute constant $c$, implying that for these $p$ the $L^p$-Brunn-Minkowski inequality holds for sufficiently close pairs of bodies \cite[Theorem 1.1]{KM}. They also treated the case of $\ell^n_q$-unit balls $B^n_q$ (which are not $C^2_+$ for $q > 2$), and proved suitably modified versions of Conjecture \ref{local_pBM_KM} for $p = 0$ and $q \in [2, \infty]$, enabling them to prove the log-Brunn-Minkowski inequality for pairs of bodies sufficiently close to $B^n_q$ \cite[Theorem 1.2, Theorem 1.3]{KM}, generalizing the results of \cite{CLM, CL}.


Using PDE methods, Chen, Huang, Li, and Liu \cite{CHLL} have extended the local results of \cite{KM} to global results, confirming Conjecture \ref{pBM_ineq} for $p \in [1 - c n^{-\frac{3}{2}}, 1]$. We explain their approach in subsection \ref{rel_work}.

\subsection{Our results}
First of all, we state a slightly different form of Conjecture \ref{local_pBM_KM}, which is more convenient for our purposes.

\begin{conjecture} \label{local_pBM} Let $p \in [0, 1)$. For any two centrally symmetric convex bodies $K, L \in \mathcal K^n$, we have
\begin{equation}\label{local_pBM_eq}
\vol(K) \cdot \left( n(n - 1) V(L[2], K[n - 2]) + (1 - p)\int \frac{h_L^2}{h_K} dS_K\right) - n(n - p) V(L, K[n - 1])^2 \le 0
\end{equation}
\end{conjecture}

The equivalence of Conjecture \ref{local_pBM} to Conjecture \ref{local_pBM_KM} follows by a standard approximation argument, which we provide, for completeness, in Appendix \ref{KM_conn}.

Our main result is the following:
\begin{theorem}\label{main_thm} For any $p \in [0, 1)$, Conjecture \ref{local_pBM} in dimension $n$ is equivalent to the $L^p$-Brunn-Minkowski inequality (Conjecture \ref{pBM_conj}) in dimension $n$.
\end{theorem}

As a corollary of this result and the work of \cite{KM}, we obtain some cases of the $L^p$-Brunn-Minkowski inequality:

\begin{corollary}Conjecture \ref{pBM_conj} holds for $p \in [1 - c n^{-\frac{3}{2}}, 1]$.
\end{corollary}

As we stated above, this result was obtained in \cite{CHLL} by a different approach, using PDE methods; for details see the next subsection.

Section \ref{pBM_equiv} is devoted to the proof of Theorem \ref{main_thm}. The proof idea is similar to that of \cite{CLM} and \cite{KM}, except that we work in the setting of strongly isomorphic polytopes rather than smooth and strongly convex bodies. The simpler behavior of these bodies enables us to go from the local inequality to the global inequality. 

Section \ref{plane_pf} proves Conjecture \ref{local_pBM} for $n = 2$, which along with Theorem \ref{main_thm} yields a new proof of the log-Brunn-Minkowski inequality in the plane.

\subsection{Related work}\label{rel_work}
To explain the connection of this work to known results, we must introduce a new ingredient (which shall play no role outside this section). The $L^p$-Brunn-Minkowski conjecture is intimately related to the uniqueness side of the $L^p$-Minkowski problem posed by Lutwak \cite{L2}, which asks about the existence and uniqueness of convex bodies $K$ with given $L^p$-surface area measure:
\begin{equation}\label{lp_mink} 
h_K^{1 - p} \,dS_K = \mu
\end{equation} 
Both sides of the problem have been studied by many authors for various ranges of $p$, both in the centrally symmetric and non-centrally-symmetric cases; see \cite{CLZ} for an updated bibliography.

Kolesnikov and Milman \cite[\S 11]{KM} showed that a slight strengthening of inequality \eqref{local_pBM_KM_eq} for a given $p \in [0, 1)$, in which it is assumed that strict inequality holds unless $z$ is a constant, implies a local uniqueness result for the $L^p$-Minkowski problem for even measures with the same $p$. They also showed that if \eqref{local_pBM_KM_eq} for a given $p_0$, then strict inequality holds in \eqref{local_pBM_KM_eq} for nonconstant $z$ and any $p > p_0$ so their work implies local uniqueness for the $L^p$-Minkowski problem for any $p \in (1 - \frac{c}{n^{\frac{3}{2}}}, 1]$.

Using PDE methods, Chen, Huang, Li, and Liu \cite{CHLL} have recently shown that global uniqueness for the $L^p$-Minkowski problem for $C^2_+$ bodies follows from the local uniqueness property defined by \cite{KM}. In addition, they showed that uniqueness for the $L^p$-Minkowski problem in the case of $C^2_+$ bodies implies the $L^p$-Brunn-Minkowski conjecture (this was proven in \cite{BLYZ} in the case $p = 0$ in the setting of general convex bodies). Hence, the result of \cite{CHLL} can also be used to obtain Theorem \ref{main_thm}: if the local $L^p$-Brunn-Minkowski conjecture holds for a given $p$, then by the work of \cite{KM} local $L^q$-Minkowski uniqueness holds for all $q \in (p, 1)$, and by \cite{CHLL}, global $L^q$-Minkowski uniqueness, and hence the global $L^q$-Brunn-Minkowski conjecture, holds for all $q \in (p, 1)$; the $L^p$-Brunn-Minkowski conjecture for the original $p$ then follows by taking the limit.

Thus, in short, while we prove a local-to-global result for the $L^p$-Brunn-Minkowski conjecture, from which uniqueness in the $L^p$-Minkowski problem may be deduced, \cite{CHLL} prove a local-to-global result for uniqueness in the $L^p$-Minkowski problem, from which the $L^p$-Brunn-Minkowski conjecture may be deduced. In addition, our methods are purely convex-geometric.

\section{Preliminaries}

We collect here the notation and basic facts in convex geometry we shall use. A comprehensive and up-to-date reference on the theory of convex bodies is the book of Schneider \cite{S}.

A convex body $K \subset \mathbb R^n$ is a compact convex set with nonempty interior. $K$ is said to be centrally symmetric if $K = -K = \{-x: x \in K\}$. We write $\mathcal K^n$ for the set of convex bodies in $\mathbb R^n$, and $\mathcal K^n_s \subset \mathcal K^n$ for the set of centrally symmetric convex bodies.

The support function $h_K : \mathbb R^n \to \mathbb R$ associated with the convex body $K$ is defined, for $u \in \mathbb R^n$, by
	\begin{equation}\label{supp}
		h_K(u) = \max\{\langle u, y\rangle : y\in K \}
	\end{equation}
The support function is convex and positively homogeneous of degree one, so it is completely determined by its restriction to the unit sphere $S^{n - 1}$.


By a classical theorem of Minkowski (see \cite[Theorem 5.1.7]{S}), for two convex bodies $K, L \subset \mathbb R^n$, the volume of $K + t L$ can be expressed as a polynomial of degree $n$ in $t$:
\begin{equation}
\vol(K + t L) = \sum_{i = 0}^n \binom{n}{i} V(K[i], L[n - i]) t^{n - i}
\end{equation}
The coefficients $V(K[i], L[n - i])$ are known as the mixed volumes of $K$ and $L$; we have $V(K[n], L[0]) = \vol(K)$, $V(K[0], L[n]) = \vol(L)$.

Let $\mathcal H^{n - 1}$ denote the $(n - 1)$-dimensional Hausdorff measure on $\mathbb R^n$. For $\mathcal H^{n - 1}$-almost every $x \in \partial K$, there exists a unique normal vector to $K$ at $x$, namely, $u \in S^{n - 1}$ such that $h_K(u) = \langle x, u\rangle$; denote this vector by $\nu_K(x)$. Thus we have an almost-everywhere-defined map $\nu_K: \partial K \to S^{n - 1}$, called the Gauss map. The surface area measure $S_K$ of a convex body $K \subset \mathbb R^n$ is a Borel measure on $S^{n - 1}$ defined by $S_K(\omega) = \mathcal H^{n - 1}(\nu_K^{-1}(\omega))$. We have \cite[Theorem 5.1.7]{S}
\begin{equation}\label{surf_vol}
\vol(K) = \frac{1}{n}\int h_K \,dS_K \qquad\text{and}\qquad V(L, K[n - 1]) = \frac{1}{n}\int h_L\,dS_K
\end{equation} 

For a continuous function $h: S^{n - 1} \to (0, \infty)$, the Wulff shape or Alexandrov body of $h$ is defined as $A[h] = \bigcap_{u \in S^{n - 1}} \{x \in \mathbb R^n: \langle x, u\rangle \le h(u)\}$. If $h$ is the support function of a convex body $K$, then $A[h]$ is simply $K$, but in general, if $K = A[h]$, we only have $h_K \le h$.

We shall use a version of Alexandrov's lemma for Wulff shapes \cite[Lemma 2.1]{BLYZ}:

\begin{lemma}\label{alex_lem} Suppose $k(t, u): I \times S^{n - 1} \to (0, \infty)$ is continuous and differentiable in the first variable, where $I \subset \mathbb R$ is an open interval. Suppose also that for $t \in I$, the convergence $\lim_{s \to 0} \frac{k(t + s, u) - k(t, u)}{s} = \frac{\partial k(t, u)}{\partial t}$ is uniform on $S^{n - 1}$. If $K_t = A[k(t,\cdot)]$ for $t \in I$, we have
\begin{equation}\label{wulff_diff}
\frac{d}{dt} \vol(K_t) = \int_{S^{n - 1}} \frac{\partial k(t, u)}{\partial t} \,dS_{K_t}
\end{equation}
In addition, the same formula holds for the one-sided derivative of $K_t$, assuming only one-sided convergence of $\frac{k(t + s, u) - k(t, u)}{s} \to \frac{\partial k(t, u)}{\partial t}$.
\end{lemma}

Let $B$ be the Euclidean ball in dimension $n$. Define the Hausdorff metric on $\mathcal K^n$ via $\delta(K, L) = \inf \{d > 0: \text{$L \subset K + d B$ and $K \subset L + dB$}\}$. The Hausdorff metric has a number of useful properties:

\begin{proposition}\label{haus_conv}\mbox{}
\begin{enumerate}
	\item A sequence of convex bodies $\{K_i\}_{i = 1}^\infty \subset \mathcal K^n$ converges to a convex body $K$ if and only if $h_{K_i}$ converges to $h_K$ uniformly on $S^{n - 1}$ \cite[Lemma 1.8.14]{S}.
	\item The mixed volumes $V(K[i], L[n - i])$ are continuous with respect to $\delta$ \cite[p. 280]{S}.
	\item The surface area measures $S_K$ are weakly continuous with respect to $\delta$ \cite[Theorem 4.2.1]{S}.
	\item If $h_i: S^{n - 1} \to (0, \infty)$ converge uniformly to $h$, then the Wulff shapes $A[h_i]$ converge in the Hausdorff metric to $A[h]$ \cite[Theorem 7.5.2]{S}.
\end{enumerate}
\end{proposition}

Finally, given two convex bodies $K, L \subset \mathbb R^n$ such that $0 \in \intr K, \intr L$, $p \in [0, 1)$ and $\lambda \in [0, 1]$, we shall write $(1 - \lambda) h_K +_p \lambda h_L$ for $((1 - \lambda) h_K^p + \lambda h_L^p)^{\frac{1}{p}}$ or $h_K^{1 - \lambda} h_L^\lambda$ in the cases $p > 0$, $p = 0$, respectively. The $L^p$-Minkowski combination of $K$ and $L$, $(1 - \lambda) K +_p L$, is defined as $A[(1 - \lambda) h_K  +_p \lambda h_L]$.

\subsection{Strongly isomorphic polytopes}\label{si_polytopes}
Our reference for the theory of polytopes is \cite[\S 2.4]{S}.

A convex polytope $P \subset \mathbb R^n$ is a convex body which can be written as the intersection of a finite set of half-spaces: there exist $u_1,\ldots,u_N \in S^{n - 1}$, $h_1(P),\ldots,h_N(P) \in \mathbb R$ such that
\begin{equation} P = \bigcap_{i = 1}^N \{x: \langle x, u_i\rangle \le h_i(P)\}
\end{equation}
We set $F_i(P) = \{x \in P: \langle x,u_i\rangle = h_i(P)\}$, the ith facet of $P$, a polytope of dimension at most $n - 1$; when there is no chance of confusion we simply write $h_i$, $F_i$. (Note that depending on the $h_i$, some of the faces may be empty.) The volume of $P$ may be computed as $|P| = \frac{1}{n} \sum_{i = 1}^N h_i |F_i|$, and the surface area measure of $P$ is a discrete measure: 
\begin{equation}\label{surf_poly} 
S_P = \sum_{i = 1}^N |F_i| \delta_{u_i}
\end{equation}

For a given facet $F_i$ of $P$ and any $j \in \{1,\ldots,N\}\backslash\{i\}$ such that $F_{ij}(P) := F_i \cap F_j$ is nonempty, we write $u_{ij} = \frac{u_j - \langle u_i, u_j\rangle u_i}{\sqrt{1 - \langle u_i, u_j\rangle^2}}$ for the normal to $F_{ij}$ relative to $F_i$; then we have 
\begin{equation}\label{u_ij} u_{j} = u_i \cos \theta_{ij} + u_{ij} \sin \theta_{ij}
\end{equation} 
where $\theta_{ij}$ is the angle between $u_i$ and $u_j$. Taking the inner product with an arbitrary $x \in F_i$, we obtain 
\begin{equation}\label{h_ij}
h_{ij}(P) := h_{F_i(P)}(u_{ij}) = h_j(P) \csc \theta_{ij} - h_i(P) \cot \theta_{ij}
\end{equation}
Note that in terms of the $h_{ij}$, we may write $|F_i| = \frac{1}{n - 1} \sum_{j = 1}^N h_{ij} |F_{ij}|$.

Two polytopes $P, Q$ are called strongly isomorphic if $\dim F(P, u) = \dim F(Q, u)$ for all $u \in S^{n - 1}$. In particular, this implies that the normal vectors $u_1,\ldots,u_N$ to $P$ and $Q$ are the same, that $F_i(P)$ is nonempty precisely when $F_i(Q)$ is, and the same for $F_{ij}(P)$, $F_{ij}(Q)$. Strong isomorphism is clearly an equivalence relation on the set of polytopes in $\mathcal K^n$; the equivalence classes under this relation are called a-types.

It follows from the proof of \cite[Theorem 5.1.7]{S} that for strongly isomorphic polytopes $P, Q$ we have 
\begin{align}
V(Q, P[n - 1]) &= \frac{1}{n}\sum_{i = 1}^N h_i(Q) |F_i(P)|  \label{mixed_vols1}
\\
V(Q[2], P[n - 2]) &= \frac{1}{n(n - 1)} \sum_{i,j = 1}^N h_i(Q) h_j(Q) \Gamma_{ij}(P) \label{mixed_vols2}
\end{align} 
where
\begin{equation}\label{gamma_ij} 
\Gamma_{ij}(P) = 
\begin{cases} -\sum_{k: F_{ik}(P) \neq\emptyset} \cot \theta_{ik} |F_{ik}| &\qquad i = j \\
               \csc \theta_{ij} |F_{ij}|                                   &\qquad i \neq j, \text{$F_{ij}(P)$ is nonempty} \\
               0                                                           &\qquad \text{otherwise}
\end{cases}
\end{equation}


A polytope $P$ is said to be simple if each vertex of $P$ is the intersection of precisely $n$ facets. We shall use the fact that any finite set of convex bodies may be approximated simultaneously in the Hausdorff metric by strongly isomorphic simple convex polytopes \cite[Theorem 2.4.15]{S}.

Given an a-type of polytopes defined by normal vectors $u_1,\ldots,u_N$ (as well as lower-dimensional intersection data), the $N$-tuple of support functions $(h_P(u_1),\ldots,h_P(u_N))$ of a polytope $P$ in the a-type is called the support vector of $P$. We have the following theorem \cite[Lemma 2.4.13]{S}:

\begin{theorem}\label{support_ngbd} Let $P$ be a simple polytope with facet normals $u_1,\ldots,u_N$ and support vector $h = (h_1,\ldots,h_N) \in \mathbb R^N$. There exists a neighborhood $U$ of $h$ such that any $h' \in U$ is the support vector of a polytope $P'$ strongly isomorphic to $P$.
\end{theorem}

\section{From the \texorpdfstring{$L^p$}{Lp}-Brunn-Minkowski inequality to its local version and back}\label{pBM_equiv}
First, we fix some notations we shall use throughout the section. Let $K, L \in \mathcal K^n_s$, and fix $p \in [0, 1]$. Write $K_\lambda = (1 - \lambda)K +_p \lambda L$ and
\begin{equation}\label{V_KL}
V_{K,L}(\lambda) = \begin{cases} \vol(K_\lambda)^{\frac{p}{n}} & p > 0 \\
             \log \vol(K_\lambda) &  p = 0
\end{cases}
\end{equation}
We recall a lemma from \cite{BLYZ}. (The lemma was proven there for $p > 0$; the proof for $p = 0$ is identical.)

\begin{lemma}\label{log_mink_conc}
The $L^p$-Brunn-Minkowski inequality (Conjecture \ref{pBM_conj}) holds if and only if for all $K, L \in \mathcal K^n_s$, $V_{K, L}(\lambda)$ is concave on $[0, 1]$.
\end{lemma}

We restate here our main conjecture, the local $L^p$-Brunn-Minkowski inequality, for the reader's convenience:

\begin{conjecture} \label{local_pBM2} For any two bodies $K, L \in \mathcal K^n_s$, we have
\begin{equation}\label{local_pBM_eq2}
\vol(K) \cdot \left( n(n - 1) V(L[2], K[n - 2]) + (1 - p)\int \frac{h_L^2}{h_K} dS_K\right) - n(n - p) V(L, K[n - 1])^2 \le 0
\end{equation}
\end{conjecture}

Our first observation is that it is sufficient to consider the case of simple strongly isomorphic polytopes in Conjecture \ref{local_pBM2}. This follows immediately by approximating $K, L$ by sequences of simple strongly isomorphic polytopes converging to $K, L$ and using the continuity of all the relevant quantities with respect to the Hausdorff metric (Proposition \ref{haus_conv}).

We shall prove that Conjecture \ref{local_pBM2} is equivalent to the $L^p$-Brunn-Minkowski inequality \eqref{pBM_ineq} in two main steps. First, in Proposition \ref{second_der}, we compute the second derivative of $V_{K, L}(\lambda)$ for strongly isomorphic polytopes $K, L$ at $\lambda \in [0, 1]$ such that the bodies $K_{\lambda'}$ for $\lambda'$ in a neighborhood of $\lambda$ are strongly isomorphic to one another, and show that in this case, $V_{K,L}''(\lambda) \le 0$ is equivalent to inequality \eqref{local_pBM_eq2} (for a different pair of bodies). The computation is similar in spirit to work of \cite{CLM} and \cite{KM}, but is technically easier in the polytope setting. Using Lemma \ref{log_mink_conc}, this shows that the $L^p$-Brunn-Minkowski inequality implies Conjecture \ref{local_pBM2}.

Next, in Proposition \ref{poly_logconc}, we argue that after removing a finite number of points, we can divide $[0, 1]$ into a disjoint union of open intervals $I_j$ such that the $K_\lambda$, $\lambda \in I_j$ are all strongly isomorphic to one another. Hence, assuming Conjecture \ref{local_pBM2} and using the computation of Proposition \ref{second_der}, we obtain that $V_{K,L}$ has nonpositive second derivative except at a finite number of points, which suffices to show that $V_{K, L}$ is concave because $V_{K, L}$ is continuously differentiable. Thus Conjecture \ref{local_pBM2} implies the log-Brunn-Minkowski inequality.

We now proceed to the details of the proof. First, we fix a few more notations. Let $K, L \in \mathcal K^n_s$ be a fixed pair of strongly isomorphic polytopes, and let $u_1,\ldots,u_N \in S^{n - 1}$ be the facet normals to $K, L$. 

\begin{lemma} For any $p \in [0, 1]$ and any $\lambda \in [0, 1]$,
\begin{equation}\label{K_lambda} K_\lambda = \{x: \langle x, u_i\rangle \le ((1 - \lambda) h_K +_p \lambda h_L)(u_i),\,i = 1,\ldots,N\}
\end{equation}
\begin{proof}
This is not immediately obvious, because $K_\lambda$ is defined to be the intersection of the half-spaces $H_{K_\lambda, v}^- = \{x: \langle x, v \rangle \le ((1 - \lambda) h_K +_p \lambda h_L)(v)\}$ for all $v \in S^{n - 1}$, and a priori, we do not know that for $v \neq u_1,\ldots,u_N$, $H_{K_\lambda, v}^-$ doesn't contribute to bounding $K_\lambda$. It is thus necessary to demonstrate that the set on the RHS of \eqref{K_lambda}, which we temporarily denote $K_\lambda'$, is already contained in $H_{K_\lambda, v}^-$, which is to say that $h_{K_\lambda'}(v) \le ((1 - \lambda) h_K +_p \lambda h_L)(v)$ for all $v \in S^{n - 1}\backslash\{u_1,\ldots,u_N\}$.

Since $v$ is not a normal vector of $K$, $v$ lies in the normal cone of a face $F$ of $K$ of codimension at least $2$, so we may write $v = \sum_{k = 1}^m c_k u_{i_k}$ for $2 \le m \le n$, $i_1,\ldots,i_m \in \{1,\ldots,N\}$, and $c_j > 0$. Hence $h_K(v) = \sum_{k = 1}^m c_j h_K(u_{i_k})$, and similarly $h_L(v)$, $h_{K_\lambda'}(v)$. So we need to prove that 
\begin{equation}\sum c_k ((1 - \lambda) h_K(u_{i_k}) +_p \lambda h_L(u_{i_k}) \le (1 - \lambda) h_K(v) +_p \lambda h_L(v)
\end{equation}

For $p > 0$, the inequality we wish to prove, after raising both sides to the power $p$, becomes
\begin{equation}\label{hkl_ineq}
\left(\sum c_k ((1 - \lambda) h_K(u_{i_k})^p + \lambda h_L(u_{i_k})^p)^{\frac{1}{p}}\right)^p \le (1 - \lambda)\left(\sum c_k h_K(u_{i_k})\right)^p + \lambda\left(\sum c_k h_L(u_{i_k})\right)^p
\end{equation}
Let $U = \{u_{i_1},\ldots,u_{i_m}\}$, and define a measure $\mu$ on $U$ via $\mu(\{u_{i_k}\}) = c_k$. Considering $h_K, h_L$ as functions on $U$, we see that \eqref{hkl_ineq} is precisely the triangle inequality $\|f + g\| \le \|f\| + \|g\|$ in $L^{\frac{1}{p}}(\mu)$ applied to  $f = (1 - \lambda) h_K^p$ and $g = \lambda h_L^p$.

The case $p = 0$ may be obtained by a limiting argument, but we also give a simple direct proof. Writing $x_k = \frac{h_L(u_{i_k})}{h_K(u_{i_k})}$, $a_j = c_k h_K(u_{i_k})$, we wish to prove 
\begin{equation}\sum a_k x_k^\lambda \le \left(\sum a_k x_k\right)^\lambda \left(\sum a_k\right)^{1 - \lambda}
\end{equation} 
By homogeneity, we may assume $\sum a_k = 1$, so it suffices to show 
\begin{equation}\label{jensen_log}
\sum a_k x_k^\lambda \le \left(\sum a_k x_k\right)^\lambda
\end{equation}
But $x \mapsto x^\lambda$ is concave for $\lambda \in [0, 1]$, so \eqref{jensen_log} follows from Jensen's inequality. Rewinding, we obtain 
\begin{equation} 
\sum c_k h_K(u_{i_k})^\lambda h_L(u_{i_k})^{1 - \lambda} \le h_K(v)^\lambda h_L(v)^{1 - \lambda}
\end{equation} 
as desired. This concludes the proof of the lemma.
\end{proof}
\end{lemma}

\begin{corollary}Let $\sigma, \tau, \mu \in [0, 1]$ and $\lambda = (1 - \mu) \sigma + \mu \tau$. If $K_\sigma, K_\lambda, K_\tau$ are all strongly isomorphic then $K_\lambda = (1 - \mu) K_\sigma +_p \mu K_\tau$.
\begin{proof}
Let $v_j$, $j = 1,\ldots,M$ be the facet normals of $K_\sigma, K_\tau$; by the lemma, $\{v_j\}_{j = 1}^M \subset \{u_i\}_{i = 1}^N$. Since $K_\lambda$ is strongly isomorphic to $K_\sigma, K_\tau$, in particular it has the same facet normals, so we may write
\begin{equation}
K_\lambda = \{x: \langle x, v_j\rangle \le ((1 - \lambda) h_K +_p \lambda h_L)(v_j),\,j = 1,\ldots,M\}
\end{equation}
Also, for every $v_j$ which is a facet normal of $K_\sigma, K_\tau$ we have $h_{K_\sigma}(v_j) = ((1 - \sigma) h_K +_p \sigma h_L)(v_j)$ and the same for $h_{K_\tau}(v_j)$, again by the lemma. Thus 
\begin{align}((1 - \mu) h_{K_\sigma} +_p \mu h_{K_\tau})(v_j) &= ((1 - \mu) ((1 - \sigma) h_K +_p \sigma h_L) +_p (\mu (1 - \tau) h_K +_p \tau h_L))(v_j) \nonumber \\
&= ((1 - \lambda) h_K +_p \lambda h_L)(v_j) 
\end{align}
for all $j$, giving $K_\lambda = (1 - \mu) K_\sigma +_p \mu K_\tau$ as desired.
\end{proof}
\end{corollary}

\begin{proposition}\label{second_der} The local $L^p$-Brunn-Minkowski conjecture (Conjecture \ref{local_pBM2}) is equivalent to the following statement:
\begin{equation}\label{second_der_eq}
\parbox[t]{0.9\textwidth}
{
  \begin{spacing}{1.25}
  {\normalfont
  For any two strongly isomorphic polytopes $K, L$ as above, and for any $\lambda \in [0, 1]$ such that there exists a (possibly one-sided) neighborhood $U$ of $\lambda$ for which all the $\{K_{\lambda'}: \lambda' \in U\}$ are strongly isomorphic to one another, we have $V_{K, L}''(\lambda) \le 0$.
  }
  \end{spacing}
}
\tag{$*$}
\end{equation}

\begin{remark}We shall see in the course of the proof that $V_{K, L}$ is twice differentiable at any $\lambda$ satisfying the conditions of \eqref{second_der_eq}.
\end{remark}

\begin{proof}First, we prove that Conjecture \ref{local_pBM2} implies \eqref{second_der_eq}. Given $\lambda, U$ satisfying the conditions of \eqref{second_der_eq}, let $[\sigma, \tau] \subset U$. Then $K_\sigma, K_\tau$ are strongly isomorphic and for any $\alpha = (1 - \mu) \sigma + \mu \tau$, $\mu \in [0, 1]$, we have $K_\alpha = (1 - \mu) K_\sigma +_p \mu K_\tau$ by the preceding corollary. Thus, we reduce to the case where all the $\{K_{\lambda'}: \lambda' \in U\}$ are strongly isomorphic to $K$ (and to $L$).

Let $h_i = h_K(u_i)$. For $p \in (0, 1)$, we write $h_L(u_i) = h_i (1 + p s_i)^{\frac{1}{p}}$, so that $((1 - \lambda) h_K +_p \lambda h_L)(u_i) = h_i (1 + \lambda p s_i)^{\frac{1}{p}}$. For $p = 0$, write $h_L(u_i) = h_i e^{s_i}$, so that $((1 - \lambda) h_K +_o \lambda h_L) = h_i e^{\lambda s_i}$. It will be convenient to write, for $p > 0$, 
\begin{align}a^{(p)}_i(\lambda) &= (1 + \lambda p s_i)^{\frac{1}{p}} \\
b^{(p)}_i(\lambda) &= (1 + \lambda p s_i)^{\frac{1 - p}{p}} \\
c^{(p)}_i(\lambda) &= (1 + \lambda p s_i)^{\frac{1 - 2 p}{p}}
\end{align} 
and for $p = 0$, $a^{(p)}_i(\lambda) = b^{(p)}_i(\lambda) = c^{(p)}_i(\lambda) = e^{\lambda s_i}$. For any $p$, we have
\begin{align} 
(a^{(p)}_i)' &= s_i b^{(p)}_i \\ 
(b^{(p)}_i)' &= (1 - p) s_i c^{(p)}_i \\ 
(b^{(p)}_i)^2 &= a^{(p)}_i c^{(p)}_i
\end{align}

In the sequel we shall suppress the dependence of $a^{(p)}_i, b^{(p)}_i, c^{(p)}_i$ on $p$. 

For any $\lambda \in U$, $K_\lambda$ is the polytope defined by normal vectors $u_1,\ldots,u_N$ and support numbers $h_{K_\lambda}(u_i) = h_i a_i(\lambda)$. By assumption, $K_\lambda$ is strongly isomorphic to $K$, so it has facets corresponding to each $u_i$; thus, its surface area measure $S_{K_\lambda}$ is given by $\sum_{i = 1}^N |F_i(K_\lambda)| \delta_{u_i}$. We may thus compute the first derivative of $\vol(K_\lambda)$ at any $\lambda \in U$ by applying Lemma \ref{alex_lem}:
\begin{align}
 \frac{d}{d\lambda}\vol(K_\lambda) &= \int_{S^{n - 1}} \left(\frac{d}{d\lambda} \left((1 - \lambda) h_K +_p \lambda h_L\right)\right) \,dS_{K_\lambda} \nonumber \\
&= 
\sum_{i = 1}^N |F_i(K_\lambda)| \frac{d}{d\lambda} (h_i a_i) = \sum_{i = 1}^N h_i s_i b_i |F_i(K_\lambda)| \label{vol_der1} 
\end{align}

In order to compute the second derivative of $\vol(K_\lambda)$, we need to know the derivative of $\vol(F_i(K_\lambda))$. But for each $i$, $F_i(K_\lambda)$ is itself a polytope strongly isomorphic to $F_i(K)$, with normal vectors $u_{ij}$ as defined by equation \eqref{u_ij}, and support numbers $h_{ij}(K_\lambda) = h_j(K_\lambda) \csc \theta_{ij} - h_i(K_\lambda) \cot \theta_{ij}$. Hence, by the same computation, we have 
\begin{equation}\label{facevol_der}
\frac{d\vol(F_i(K_\lambda))}{d\lambda} = \sum_{j: F_{ij} \neq \emptyset} (s_j h_j  \csc \theta_{ij} b_j - s_i h_i \cot \theta_{ij}  b_i) |F_{ij}(K_\lambda)|
\end{equation}

Differentiating \eqref{vol_der1} and using \eqref{facevol_der}, we obtain 
\begin{equation}\label{vol_der2} \frac{d^2 \vol(K_\lambda)}{d\lambda^2} = (1 - p)\sum_{i = 1}^N s_i^2 h_i c_i |F_i(K_\lambda)| + \sum_{i,j = 1}^N h_i h_j s_i s_j b_i b_j \Gamma_{ij}(K_\lambda)
\end{equation}
where the $\Gamma_{ij}$ are defined by equation \eqref{gamma_ij}.

Recalling the definition of $V_{K, L}$ (Equation \eqref{V_KL}), for any $p \in [0, 1]$, we have
\begin{equation}\label{p_conc} n(n-p) \vol(K_\lambda) V_{K, L}''(\lambda) = n \vol(K_\lambda)'' \vol(K_\lambda) - (n - p) (\vol(K_\lambda)')^2
\end{equation}
Examining equations \eqref{vol_der1} and \eqref{vol_der2}, we see that the RHS may be written as $\Psi(Z) \le 0$, where $Z = (h_1 s_1 b_1(\lambda),\ldots,h_N s_N b_N(\lambda)) \in \mathbb R^N$ and $\Psi$ is the quadratic form defined by
\begin{multline}\label{psi_def}
\Psi(X) = n\left((1 - p)\sum_{i = 1}^N \frac{X_i^2}{h_i a_i} |F_i(K_\lambda)| X_i Y_i + \sum_{i, j = 1}^N \Gamma_{ij}(K_\lambda) X_i X_j \right) \vol(K_\lambda) \\ - (n - p)\left(\sum_{i = 1}^N X_i |F_i(K_\lambda)| \right)^2
\end{multline}

Let $M$ be the matrix associated to $\Psi$, and let $X_{K_\lambda} = (h_1 a_1(\lambda),\ldots,h_N a_N(\lambda))$, the support vector of the polytope $K_\lambda$. We claim that $M X_{K_\lambda} = 0$. Indeed, we have 
\begin{multline}(M X_{K_\lambda})_i = n\left((1 - p) |F_i(K_\lambda)| + \sum_{j = 1}^N h_j a_j \Gamma_{ij}(K_\lambda)\right) \vol(K_\lambda) \\ - (n - p) |F_i(K_\lambda)| \sum_{j = 1}^N h_j a_j |F_j(K_\lambda)|  
\end{multline}

We recognize the second sum as $n \vol(K_\lambda)$ and the first sum as $(n - 1) |F_i(K_\lambda)|$, so the entire expression cancels out, as desired.

Now, for $X_P$ which is the support vector of a polytope $P$ strongly isomorphic to $K$, combining formulas \eqref{surf_poly}, \eqref{mixed_vols1}, \eqref{mixed_vols2} shows that
\begin{multline}\Psi(X_P) = n\vol(K_\lambda) \left( (n(n - 1) V(P[2], K_\lambda[n - 2]) + (1 - p)\int \frac{h_P^2}{h_{K_\lambda}} dS_{K_\lambda}\right) \\ - n^2 (n - p) V(P, K_\lambda[n - 1])^2
\end{multline}

Hence, assuming Conjecture \ref{local_pBM2} yields that $\Psi(X_P) \le 0$ for any $P$ strongly isomorphic to $K$. Since, by Theorem \ref{support_ngbd}, any vector in a neighborhood of $X_{K_\lambda}$ is the support vector of a polytope strongly isomorphic to $K$, and as we computed above, $\Psi(X_{K_\lambda} + X) = \Psi(X)$ for any $X$, we see that in fact $\Psi$ is a negative semidefinite quadratic form, and hence that $V_{K, L}''(\lambda) = \Psi(X_{K_\lambda}) \le 0$. This proves one direction of the proposition.

Conversely, we wish to prove that \eqref{local_pBM_eq2} holds for two given polytopes $P, Q$ assuming \eqref{second_der_eq}. Let $u_1,\ldots,u_N$ be the normal vectors and $X_P = (h_1,\ldots,h_N)$ and $X_Q = (h'_1,\ldots,h'_N)$ the support vectors of $P$ and $Q$, respectively. Now set $K = P$, and let $L$ be the polytope strongly isomorphic to $P$ defined by the vector 
\begin{equation} 
v_\epsilon = \begin{cases} \left(h_1 \left(1 + p\epsilon\frac{h'_1}{h_1}\right)^{\frac{1}{p}}, \ldots, h_N \left(1 + p\epsilon\frac{h'_N}{h_N}\right)^{\frac{1}{p}}\right) \qquad & p > 0 \\
\left(h_1 e^{\epsilon\frac{h'_1}{h_1}}, \ldots, h_N e^{\epsilon\frac{h'_N}{h_N}}\right) & p = 0
\end{cases}
\end{equation}

By Theorem \ref{support_ngbd}, by taking $\epsilon$ sufficiently small we obtain that $v_\epsilon$ is the support vector of a polytope strongly isomorphic to $K$; moreover, $K_\lambda = (1 - \lambda) K +_p \lambda L$ is defined by the vector $v_{\lambda\epsilon}$ and hence is strongly isomorphic to $K$ for all $\lambda \in [0, 1]$. Thus we can use the computation above to obtain $V_{K, L}''$: for $L$ defined as above, we obtain $s_i = \epsilon\frac{h'_i}{h_i}$, and applying \eqref{psi_def} shows that
\begin{equation}V_{K, L}''(0) = \Psi\left(\epsilon \frac{h'_1}{h_1} \cdot h_1, \ldots, \epsilon \frac{h'_N}{h_N} \cdot h_N\right) = \epsilon^2\Psi(X_Q)
\end{equation}
But we have 
\begin{equation}
\Psi(X_Q) = \vol(P) \left( (n(n - 1) V(Q[2], P[n - 2]) + \int \frac{h_Q^2}{h_{P}} dS_{P}\right) - n^2 V(Q, P[n- 1])^2
\end{equation}
and since $V_{K, L}''(0) \le 0$ by assumption, the proof is complete.
\end{proof}
\end{proposition}

By Lemma \ref{log_mink_conc}, assuming the $L^p$-Brunn-Minkowski inequality, we have that for any two strongly isomorphic polytopes $K, L$, $V_{K, L}''(\lambda) \le 0$ everywhere it is defined, and in particular at any $\lambda$ satisfying the conditions of \eqref{second_der_eq}. Thus, by Proposition \ref{second_der}, the $L^p$-Brunn-Minkowski inequality implies Conjecture \ref{local_pBM2}.

To go in the reverse direction, from Conjecture \ref{local_pBM2}, to the $L^p$-Brunn-Minkowski inequality, requires a more careful analysis of the $L^p$-Minkowski sum of strongly isomorphic polytopes not just on well-behaved neighborhoods, but for all $\lambda \in [0, 1]$. This is the content of the following proposition.

\begin{proposition}\label{poly_logconc}
Assume that for a given $p \in [0, 1)$, Conjecture \ref{local_pBM2} holds for all $K, L \in \mathcal K^n_s$. Then for any two strongly isomorphic polytopes $K, L \in \mathcal K^n_s$, $\vol(K_\lambda)$ is $\frac{p}{n}$-concave on $[0, 1]$.
\begin{proof}By Proposition \ref{second_der}, we know that $f(\lambda) = \vol(K_\lambda)$ is $\frac{p}{n}$-concave on any subinterval $(\sigma, \tau) \subset [0, 1]$ such that all the $K_\alpha$ for $\alpha \in (\sigma, \tau)$ are strongly isomorphic. The additional ingredient we need is the following claim:

\begin{claim}There exist a finite number of disjoint open intervals $I_1, \ldots, I_m \subset [0, 1]$ such that $[0, 1] \backslash \bigcup_{j = 1}^m I_j$ is a finite set, and for each $j$, all the polytopes $K_\lambda$ for $\lambda \in I_j$ are strongly isomorphic.
\end{claim}

Deferring for the moment the proof of this claim, we show how it implies the proposition. By Lemma \ref{alex_lem}, $\vol(K_\lambda)$ is differentiable on $[0, 1]$, with derivative
\begin{equation}
\int_{S^{n - 1}} \frac{d}{d\lambda} ((1 - \lambda) h_K +_p \lambda h_L)\,dS_{K_\lambda} = 
\begin{cases} \displaystyle\int_{S^{n - 1}} \frac{h_L^p - h_K^p}{p} ((1 - \lambda) h_K^p + \lambda h_L^p)^{\frac{1 - p}{p}} \,dS_{K_\lambda} & p > 0 
\vspace{10pt} 
\\ 
\displaystyle\int_{S^{n - 1}} (h_K^{1 - \lambda} h_L^\lambda) \log \frac{h_L}{h_K} \,dS_{K_\lambda}
& p = 0
\end{cases}
\end{equation}

For any $p$, the integrand is $L^\infty$-continuous in $\lambda$ and $S_{K_\lambda}$ is weakly continuous in $\lambda$, so $\vol(K_\lambda)$ is continuously differentiable, and the same is true for $V_{K, L}(\lambda)$. Hence, using the claim, $V_{K, L}'(\lambda)$ is a continuous function which is nonincreasing outside of a finite set, so it must be nonincreasing on all of $[0, 1]$; in other words, $\vol(K_\lambda)$ is $\frac{p}{n}$-concave on $[0, 1]$, as desired.

We now turn to proving the claim above. As usual, let $u_1,\ldots,u_N$ be the facet normals of $K, L$, and let $h_i(\lambda) = h_K(u_i)^{1 - \lambda} h_L(u_i)^\lambda$ for $i = 1,\ldots,N$ be the corresponding support numbers. First of all, we ask when the facet $F_i(K_\lambda) = \{x \in K_\lambda: \langle x, u_i\rangle = h_i(K_\lambda)\}$ becomes empty. Clearly, this holds whenever 
\begin{equation}\label{facet_vanish}
h_i(\lambda) > h_{K_\lambda}(u_i)
\end{equation}
In this case, there exists a vertex $v$ of $K_\lambda$ such that $u_i$ lies in its normal cone, which means that there exist $n$ linearly independent facet normals $u_{i_1},\ldots,u_{i_n}$, $i_j \neq i$, such that $\langle v, u_{i_k}\rangle = h_{i_j}(\lambda)$, and $u_i$ can be written as a sum $\sum_{j = 1}^n c_j u_{i_k}$ where the $c_j$ are all nonnegative. In particular, we have 
\begin{equation}
h_{K_\lambda}(u_i) = \langle v, u_i\rangle = \sum_{j = 1}^n c_j h_{K_\lambda}(u_{i_k})
\end{equation}
and hence inequality \eqref{facet_vanish} may be written as
\begin{equation}
h_i(\lambda) > \sum_{j = 1}^n c_j h_{i_j}(\lambda)
\end{equation}

Note that for different $\lambda \in [0, 1]$, the vertex structure of $K_\lambda$ will in general differ. Thus, we must consider the inequality $h_i > \sum_{j = 1}^n c_j h_{i_j}$ not just for one set of $u_{i_1},\ldots,u_{i_N}$, but for any set of $n$ vectors $u_{i_1},\ldots,u_{i_N}$ such that $u_i = \sum_{j = 1}^n c_j u_{i_k}$ with nonnegative $c_j$, since at some $\lambda \in [0, 1]$, the corresponding facets $F_{i_j}$ may form a vertex $v$, in whose normal cone $u_i$ will necessarily lie. Conversely, if $h_i(\lambda) < \sum_{j = 1}^n c_j h_{i_j}(\lambda)$ for all $u_{i_1},\ldots,u_{i_N}$ in whose positive cone $u_i$ lies, then $u_i$ cannot lie in the relative interior of the normal cone of any lower-dimensional face of $K_\lambda$, implying that $u_i$ is a facet normal.  Hence, there is a finite set of equations of the form $h_i = \sum c_j h_{i_j}$ such that the facet structure of $K_\lambda$ can change only at their solutions.

The situation is made slightly more complicated by the fact that varying the support numbers of a polytope may cause lower-dimensional faces to disappear without any of the facets vanishing. However, we can think of the $(n - 2)$-dimensional faces of a given facet $F_i$ as facets of $F_i$ and apply the preceding argument. As equation \eqref{h_ij} shows, the support numbers of the $(n - 2)$-dimensional faces $F_{ij} = F_i \cap F_j$ bounding a given facet are linear combinations (with coefficients independent of $\lambda$) of the $h_i$, so the condition for vanishing of $F_{ij}$ (given that $F_i, F_j$ do not vanish) may also be expressed as a linear inequality in the $h_i$. By induction, the same holds true for all lower-dimensional faces of $K_\lambda$.

The upshot of this discussion is that there is a finite set of equations $\sum_{j = 1}^n c_{ij} h_j(\lambda) = 0$, $i = 1,\ldots,M$ such that $K_\lambda$'s a-type changes only at solutions of these equations; it thus suffices to show that each such equation has at most finitely many solutions. Each of these equations may be written as 
\begin{equation}\label{supp_sum}
\sum_{j = 1}^N c_{ij}' a^{(p)}_i(\lambda) = 0
\end{equation}
where $c_{ij}' = c_{ij} h_K(u_j)$ and $a^{(p)}_i(\lambda)$ is defined as in Proposition \ref{second_der}.  The $a^{(p)}_i(\lambda)$ are analytic on a neighborhood of $[0, 1]$, and hence so is $\sum_{j = 1}^N c_{ij}' a^{(p)}_i(\lambda)$; thus, if this expression does not vanish identically, it has a finite number of zeros on $[0, 1]$. If it does vanish, this means that the face corresponding to the equation $\sum_{j = 1}^n c_{ij} h_j(\lambda) = 0$ is degenerate (if nonempty) for all $\lambda \in [0,1]$, and hence need not be considered. This concludes the proof of the claim, and with it, the proof of Proposition \ref{poly_logconc}.
\end{proof}
\end{proposition}

Now, assume Conjecture \ref{local_pBM2} holds, and let $K, L$ be any two bodies in $\mathcal K^n_s$. We can approximate $K, L$ by sequences of strongly isomorphic polytopes $P_m, Q_m$ converging to $K, L$ respectively. For every $\lambda \in [0, 1]$, 
\begin{equation}(1 - \lambda) P_m +_p \lambda Q_m \to (1 - \lambda)K +_p \lambda L
\end{equation}
by Proposition \ref{haus_conv} (d), and hence 
\begin{equation}\vol((1 - \lambda) P_m +_p \lambda Q_m) \to \vol((1 - \lambda)K +_p \lambda L)
\end{equation} 
for all $\lambda \in [0, 1]$. Since a pointwise limit of $\frac{p}{n}$-concave functions is $\frac{p}{n}$-concave, $\vol((1 - \lambda)K +_p \lambda L)$ is $\frac{p}{n}$-concave. This completes the proof that Conjecture \ref{local_pBM2} implies the $L^p$-Brunn-Minkowski inequality, and the proof of Theorem \ref{main_thm} is complete.

\section{The proof of the log-Brunn-Minkowski inequality in the plane}\label{plane_pf}
The log-Brunn-Minkowski inequality in the plane follows immediately from Theorem \ref{main_thm} and the following theorem:

\begin{theorem} Conjecture \ref{local_pBM} is true for $p = 0$ and $n = 2$.
\begin{proof}Let $K, L \in \mathcal K^2_s$. Define
\begin{align}
r(L, K) = \sup \{t \ge 0: \text{there exists $x \in \mathbb R^n$ such that $x + t K \subset L$}\} \\
R(L, K) = \inf \{t \ge 0: \text{there exists $x \in \mathbb R^n$ such that $x + t K \supset L$}\}
\end{align}
the inradius and circumradius of $L$ with respect to $K$, respectively. Blaschke's extension of the Bonnesen inequality \cite[Lemma 4.1]{BLYZ} states that for any plane convex bodies $K, L$ and $t \in [r(L, K), R(L, K)]$, we have
\begin{equation}\label{bonn_plane}
\vol(L) - 2 t V(L, K) + t^2 \vol(K) \le 0
\end{equation}
Since $K, L$ are centrally symmetric, we have $r(L, K) = \min_{u \in S^{n - 1}} \frac{h_L(u)}{h_K(u)}$ and $R(L, K) = \max_{u \in S^{n - 1}} \frac{h_L(u)}{h_K(u)}$, and hence
\begin{equation}
\vol(L) - 2 \frac{h_L(u)}{h_K(u)} V(K, L) + \left(\frac{h_L(u)}{h_K(u)}\right)^2 \vol(K) \le 0
\end{equation}
for all $u \in S^{n - 1}$. Integrate this inequality over the measure $h_K dS_K$ to obtain

\begin{equation}
\vol(L)\cdot \int h_K\, dS_K - 2 \int h_L\,dS_K V(K, L) + \int \frac{h_L^2}{h_K}\,dS_K \cdot \vol(K) \le 0
\end{equation}

By equation \eqref{surf_vol}, this reduces to
\begin{equation}
2\vol(K)\vol(L) - 4 V(K, L)^2 + \vol(K) \int \frac{h_L^2}{h_K}\,dS_K \le 0
\end{equation}
which is precisely \eqref{local_pBM_eq} in dimension $2$ for $p = 0$.
\end{proof}
\end{theorem}

\begin{remark}B\"or\"oczky, Lutwak, Yang and Zhang \cite{BLYZ} also use Blaschke's Bonnesen-type inequality in proving the log-Brunn-Minkowski inequality in the plane, but their argument is completely different: they reduce the log-Brunn-Minkowski inequality to demonstrating that the cone-volume measure $h_K dS_K$ of a centrally symmetric convex body is unique, and prove the uniqueness of cone-volume measures for centrally symmetric convex bodies in the plane by a compactness argument, which relies crucially on an estimate obtained by integrating equation \eqref{bonn_plane} against a different measure. Our proof, in contrast, bypasses the uniqueness of cone-volume measures entirely.
\end{remark}

\appendix

\section{Equivalence of Conjecture \ref{local_pBM_KM} and Conjecture \ref{local_pBM}}\label{KM_conn}
Obtaining our formulation of the local $L^p$-Brunn-Minkowski inequality (Conjecture \ref{local_pBM}) from that of Kolesnikov and Milman (Conjecure \ref{local_pBM_KM}) is simple: for any $C^2_+$ convex body $L$, setting $z = \frac{h_L}{h_K}$ in \eqref{local_pBM_KM_eq} yields
\begin{equation} 
\frac{1}{\vol(K)} V(h_L[1], K[n - 1])^2 \geq \frac{n-1}{n-p} V(h_L[2], K[n - 2]) + \frac{1-p}{n-p} V\left(\frac{h_L^2}{h_K}[1], K[n - 1]\right)
\end{equation}
Now using the fact that the mixed volume of $C^2$ functions on the sphere coincides with the ordinary mixed volume when the functions are support functions of convex bodies, as well as the fact that $V(f[1], K[n - 1]) = \int f\,dS_K$ for any $f \in C^2(S^{n - 1})$, immediately yields \eqref{local_pBM_eq} for any $C^2_+$ and centrally symmetric $K, L$ (upon multiplying by $n(n - p) \vol(K)$). Since any body in $\mathcal K^n_s$ may be approximated by $C^2_+$ bodies in $\mathcal K^n_s$, and the quantities in \eqref{local_pBM} are continuous with respect to the Hausdorff metric, this gives Conjecture \ref{local_pBM} for any two bodies in $\mathcal K^n_s$.

To go in the converse direction requires only slightly more work. Let $K$ be a $C^2_+$ centrally symmetric convex body. By a generalization of \cite[Lemma 1.7.8]{S}, any even $f \in C^2(S^{n - 1})$ may be written as $h_L - c h_K$ for some $c > 0$ and $C^2_+$ convex body $L \in \mathcal K^n_s$. Applying this decomposition to $z h_K$, substituting in \eqref{local_pBM_KM_eq}, and using multilinearity of the mixed volume, we see that we must prove
\begin{align}
& \frac{1}{\vol(K)} V(h_L[1], K[n - 1])^2 - \frac{n-1}{n-p} V(h_L[2], K[n - 2]) - \frac{1-p}{n-p} V\left(\frac{h_L^2}{h_K}[1], K[n - 1]\right)  \\ 
  &- 2 \cdot \left(\frac{1}{\vol(K)} V(h_L[1], K[n - 1]) V(c h_K[1], K[n - 1]) - \frac{n-1}{n-p} V(h_L[1], c h_K[1], K[n - 2]) \right.  \nonumber
\\
&\qquad\qquad \left.- \frac{1-p}{n-p} V(c h_L[1], K[n - 1])\right) \nonumber  \\ 
&+ \frac{1}{\vol(K)} V(c h_K[1], K[n - 1])^2 - \frac{n-1}{n-p} V(c h_K[2], K[n - 2]) - \frac{1-p}{n-p} V(c^2 h_K[1], K[n - 1]) \ge 0 \nonumber 
\end{align}
The second and third lines are just $-2c(1 - \frac{n-1}{n-p} - \frac{1 - p}{n - p}) V(L[1], K[n - 1]) = 0$ and the fourth line is $c^2 (1 - \frac{n-1}{n-p} - \frac{1 - p}{n - p}) \vol(K) = 0$, so what remains is 
\begin{equation}
\frac{1}{\vol(K)} V(h_L[1], K[n - 1])^2 - \frac{n-1}{n-p} V(h_L[2], K[n - 2]) - \frac{1-p}{n-p} V\left(\frac{h_L^2}{h_K}[1], K[n - 1]\right) \ge 0
\end{equation}
which is \eqref{local_pBM_eq} divided by $n(n - p)\vol(K)$. Thus Conjecture \ref{local_pBM} implies Conjecture \ref{local_pBM_KM}.

\section*{Acknowledgments}
This paper is a part of the author's thesis, being carried out under the supervision of Prof. Shiri Artstein-Avidan at Tel Aviv University. I also wish to thank Shibing Chen and Emanuel Milman for helpful comments and remarks.

%

\bibliographystyle{amsplain}

\end{document}